\documentclass[11pt]{amsart}

\usepackage{amsthm,amsmath,amssymb,amscd}
\usepackage{amssymb}
\usepackage{graphics}
\newtheorem{theorem}{Theorem}[section]
\newtheorem{corollary}[theorem]{Corollary}

\newtheorem{lemma}[theorem]{Lemma}
\newtheorem{proposition}[theorem]{Proposition}
\newtheorem{definition}{Definition}[section]

\numberwithin{equation}{section}

\begin{document}
\title [Hamiltonian L-stability of Lagrangian Translating Solitons]
 {Hamiltonian L-stability of Lagrangian Translating Solitons}

\author{Liuqing Yang}

\address{Liuqing Yang, Beijing International Center for Mathematical Research, Peking University, Beijing 100871, P. R. China.}
\email{yangliuqing@math.pku.edu.cn}


\date{}

\maketitle

\begin{abstract}
In this paper, we compute the first and second variation formulas for
the F-functional of translating solitons and study the
Hamiltonian L-stability of Lagrangian translating solitons. We prove
that any Lagrangian translating soliton is Hamiltonian L-stable.
\end{abstract}

\vspace{.2in}


\allowdisplaybreaks

\vspace{.4in}

\section{Introduction}

\noindent

An $n$-dimensional submanifold $\Sigma^n$ of $\textbf{R}^{n+p}$ is
called a self-shrinker if it is the time $t=-1$ slice of a
self-shrinking mean curvature flow that disappears at $(0,0)$, i.e.
of a mean curvature flow satisfying $\Sigma_t=\sqrt{-t}\Sigma_{-1}$.
We can also consider a self-shrinker as a submanifold that satisfies
\begin{eqnarray*}
  H=-\frac{1}{2}x^\perp.
\end{eqnarray*}

An $n$-dimensional submanifold $\Sigma^n$ of $\textbf{R}^{n+p}$ is
called a translating soliton if there is a constant vector $T$ so
that $\Sigma_t=\Sigma+tT$ is a solution to the mean curvature flow. We
can also consider a translating soliton as a submanifold that
satisfies
\begin{eqnarray*}
  H=T^\perp.
\end{eqnarray*}

According to the blow up rate of the second fundamental form,
Huisken \cite{Hu2} classified the singularities of mean curvature
flows into two types: Type I and Type II. Any Type I singularity of the mean curvature
flow must be a self-shrinker (\cite{Hu2}). Type II singularity is one class of eternal solutions, which is defined for $-\infty<t<\infty$.
One of the most important example of Type II singularity is the translating soliton (\cite{Ha,HS}).

In this paper, we mainly study the stability (in some sense) of translating solitons. It was motivated by the work of Colding-Minicozzi \cite{CM}, where they introduced the concept of
F-stability of a self-shrinker in the hypersurface case. The
definitions of many concepts in their paper can be naturally
generalized to the higher codimension case (cf. \cite{ALW,AS,LL}).

Given $x_0\in \textbf{R}^{n+p}$ and $t_0>0$, $F_{x_0, t_0}$ is
defined by
\begin{eqnarray*}
  F_{x_0,t_0}(\Sigma)=(4\pi t_0)^{-\frac{n}{2}}\int_\Sigma e^{-\frac{|x-x_0|^2}{4t_0}}d\mu.
\end{eqnarray*}
In \cite{CM}, Colding-Minicozzi proved that self-shrinkers are the
critical points for the $F_{0,1}$ functional by computing the first
variation formula of $F_{0,1}$. They also computed the second
variation formula, and defined F-stability of a self-shrinker by
modding out translations. They showed that the round sphere and
hyperplanes are the only F-stable self-shrinkers in
$\textbf{R}^{n+1}$.

In 2002, Andrews-Li-Wei \cite{ALW}, Arezzo-Sun \cite{AS} and Lee-Lue
\cite{LL} independently generalized some of Colding-Minicozzi's work
\cite{CM} from the hypersurface case to the higher codimensional
case. They computed the first and second variation formulas for the
F-functional, and studied F-stability of self-shrinkers in higher
codimension.

Recently, motivated by an observation by Oh\cite{O},
Li-Zhang\cite{LZ} and the author \cite{Y} studied Lagrangian
F-stability and Hamiltonian F-stability of Lagrangian
self-shrinkers, and proved characterization theorems for Hamiltonian
F-stability of Lagrangian self-shrinkers, which characterize the
Hamiltonian F-stablity by the eigenvalues and eigenspaces of the
drifted Laplacian.

With the above known results for self-shrinkers, it is natural to
think that translating solitons might also have some similar
properties. In fact, translating solitons are also critical points for an $F$-functional, which was studied by some people (See \cite{AS2}, \cite{Sh} and \cite{Zh} for example ). The
$F$-functional is defined by
\begin{eqnarray*}
  F(\Sigma)=\int_\Sigma e^{\langle T, x\rangle} d\mu.
\end{eqnarray*}
It is not as good as in the self-shrinker case because here
$F(\Sigma)$ is usually infinity if $\Sigma$ is a translating
soliton, since any translating soliton is noncompact and $e^{\langle
T, x\rangle}\to\infty$ very quickly as $x\to\infty$. This makes it
hard to get many corresponding results as in the self-shrinker case.

However, if we require variation vector fields to have compact
support, we can still compute the first and second variation
formulas and consider stability of translating solitons. In
\cite{Sh}, Shahriyari defined L-stability of translating surfaces in $\textbf{R}^3$, and
proved that any translating graph in $\textbf{R}^3$ is L-stable.

Especially, if a translating soliton is also a Lagrangian
submanifold of the Euclidean space, we call it a Lagrangian
translating soliton. We consider Hamiltonian L-stability (see
section 2 for the definition) of Lagrangian translating solitons.
Our main theorem is
\begin{theorem}
  Any Lagrangian translating soliton is Hamiltonian
  L-stable.
\end{theorem}

Since we have this theorem, an interesting question is that whether
this theorem has some application that could help us study the
Lagrangian mean curvature flow. Besides, we are also interested in
the L-stability of translating solitons in the hypersurface case.

\vspace{.1in}

\noindent \emph{Acknowledgement: The author would like to thank
Professor William P. Minicozzi II for introducing this problem to
her and for helpful suggestions. She would also like to thank
Professor Gang Tian for his support and comments.}

\vspace{.2in}

\section{Variation formulas and Hamiltonian L-stability}
\subsection{First and second variation formulas}
Recall that the F-functional is defined by
\begin{eqnarray*}
  F(\Sigma)=\int_\Sigma e^{\langle T, x\rangle}d\mu.
\end{eqnarray*}
The first variation formula of $F$ is
\begin{lemma}
  Let $\Sigma_s\subset\textbf{R}^{n+p}$ be a compactly supported variation of $\Sigma$ with normal variation vector field $V$, then
  \begin{eqnarray}\label{F'}
    \frac{\partial}{\partial s}(F(\Sigma_s))=\int_\Sigma\langle T^\perp-H, V\rangle e^{\langle T, x\rangle}d\mu_\Sigma.
  \end{eqnarray}
\end{lemma}
{\it Proof.} From the first variation formula (for area), we know
that
\begin{eqnarray*}
  (d\mu)'=-\langle H, V\rangle d\mu.
\end{eqnarray*}
It follows that
\begin{eqnarray*}
  \frac{\partial}{\partial s}(F(\Sigma_s))=\int_\Sigma e^{\langle T,
  x\rangle}\langle T, V\rangle d\mu-\int_\Sigma e^{\langle T,
  x\rangle}\langle H, V\rangle d\mu=\int_\Sigma \langle T^\perp-H,
  V\rangle e^{\langle T, x\rangle}d\mu.
\end{eqnarray*}
This proves the lemma.  \hfill Q.E.D.

It follows that \begin{proposition}
  $\Sigma$ is a critical point for $F$ if and only if $H=T^\perp$.
\end{proposition}
The second variation formula at a critical point is
\begin{theorem}
  Suppose that $\Sigma$ is a critical point for $F$. If $\Sigma_s$ is a compactly supported normal variation of $\Sigma$, and
\begin{eqnarray*}
  \partial_s\Big|_{s=0}\Sigma_s=V,
\end{eqnarray*}
then setting $F''=\partial_{ss}\Big|_{s=0}(F(\Sigma_s))$ gives
\begin{eqnarray}\label{F''}
  F''=\int_\Sigma-\langle V, LV\rangle e^{\langle T, x\rangle}d\mu,
\end{eqnarray}
where
\begin{eqnarray*}
  LV&=&\Delta^\perp V+\nabla_{T^T}^\perp V+\big\langle \langle A, V\rangle, A\big\rangle\nonumber\\
  &=&\left(\Delta V^\alpha+\langle T,\nabla V^\alpha\rangle+g^{ik}g^{jl}V^\beta h^\beta_{ij}h^\alpha_{kl}\right)e_\alpha.
\end{eqnarray*}
\end{theorem}
{\it Proof. } Letting primes denote derivatives with respect to $s$
at $s=0$, differentiating (\ref{F'}) gives
\begin{eqnarray}\label{F2}
  F''&=&\int_\Sigma\left\{\frac{\partial}{\partial s}\Big|_{s=0}\left(\langle T-H,
  V\rangle\right)+\langle T^\perp-H, V\rangle^2\Big|_{s=0}\right\} e^{\langle T, x\rangle}\nonumber\\
  &=&\int_\Sigma\left\{-\langle H', V\rangle+\langle T-H, V'\rangle\right\} e^{\langle T, x\rangle}.
\end{eqnarray}
Similar to the derivation of the second variation formula for the
area, we have
\begin{eqnarray}\label{H'}
  \langle H', V\rangle=\langle \Delta^\perp V+g^{ik}g^{jl}V^\beta
  h^\beta_{ij}h^\alpha_{kl}e_\alpha, V\rangle.
\end{eqnarray}
On the other hand, since $\langle\left[V, T^T\right], V\rangle=0$,
it follows that
\begin{eqnarray}\label{V'}
\langle T-H, V'\rangle&=&\langle T-H, \overline\nabla_V^T
V\rangle=\langle T, \overline\nabla_V^T V\rangle=\langle T^T,
\overline\nabla_V V\rangle=-\langle \overline\nabla_VT^T,
V\rangle\nonumber\\&=&-\langle \overline\nabla_{T^T}
V,V\rangle=-\langle \nabla_{T^T}^\perp V,V\rangle.
\end{eqnarray}
Putting (\ref{H'}) and (\ref{V'}) into (\ref{F2}) gives (\ref{F''}).
This proves the theorem. \hfill Q.E.D.

\subsection{Properties of $\mathcal{L}$ and $L$}
Notice that $T=T^T+H$, where $T^T=\nabla\langle T,x\rangle$, the
linear operator defined by
\begin{eqnarray*}
  {\mathcal{L}}v=\Delta v+\langle T, \nabla v\rangle=e^{-\langle T,x\rangle}div_\Sigma\left(e^{\langle T,x\rangle}\nabla v\right)
\end{eqnarray*}
is self-adjoint in a weighted $L^2$ space. This follows immediately from Stokes' theorem. More precisely,
\begin{lemma}\label{uLv}
  If $\Sigma\subset \textbf{R}^{n+p}$ is a submanifold of $\textbf{R}^{n+p}$, $u$ is a $C^1$ function with compact support, and $v$ is a $C^2$ function, then
  \begin{eqnarray}
    \int_\Sigma u(\mathcal{L}v)e^{\langle T, x\rangle}=-\int_{\Sigma}\langle \nabla v, \nabla u\rangle e^{\langle T, x\rangle}
  \end{eqnarray}
\end{lemma}
The next corollary is an extension of Lemma \ref{uLv}, which follows
immediately by choosing cut-off functions and using the dominated
convergence theorem, the same as in the proof of Corollary 3.10 in
\cite{CM}.
\begin{corollary}
  Suppose that $\Sigma\subset \textbf{R}^{n+p}$ is a complete submanifold of $\textbf{R}^{n+p}$ without boundary. If $u$, $v$ are $C^2$ functions with
  \begin{eqnarray*}
    \int_\Sigma(|u\nabla v|+|\nabla u||\nabla v|+|u{\mathcal{L}}v|)e^{\langle T, x\rangle}<\infty,
  \end{eqnarray*}
  then we get
  \begin{eqnarray*}
    \int_\Sigma u({\mathcal{L}}v)e^{\langle T, x\rangle}=-\int_\Sigma \langle \nabla v, \nabla u\rangle e^{\langle T,x\rangle}.
  \end{eqnarray*}
\end{corollary}

Now we calculate some equalities that we think will be useful in the future.
\begin{proposition} If $\Sigma^n\subset \textbf{R}^{n+p}$ is a translating soliton, then for every constant vector field $y$,
\begin{eqnarray}\label{Ly}
  Ly^\perp=0.
\end{eqnarray}
Especially, choosing $y=T$, we have
\begin{eqnarray}\label{LH}
  LH=0.
\end{eqnarray}
Besides, we have
\begin{eqnarray}\label{LxA}
  {\mathcal{L}}x^A=T^A.
\end{eqnarray}
\end{proposition}
{\it Proof. } Fix $p\in\Sigma$ and choose an orthonormal frame
$\{e_i\}$ such that $\nabla_{e_i}e_j(p)=0$, $g_{ij}=\delta_{ij}$ in
a neighborhood of $p$. We have
\begin{eqnarray}\label{nablay}
  \nabla_{e_i}^\perp y^\perp=\nabla_{e_i}^\perp(y-\langle y, e_j\rangle e_j)=-\langle y, e_j\rangle h_{ij}^\alpha e_\alpha.
\end{eqnarray}
Especially, choosing $y=T$, we have
\begin{eqnarray}
  \nabla_{e_i}^\perp H=\nabla_{e_i}^\perp T^\perp=-\langle T, e_j\rangle h_{ij}^\alpha
  e_\alpha,
\end{eqnarray}
i.e.,
\begin{eqnarray}
  H^\alpha_{,i}=-\langle T, e_j\rangle h_{ij}^\alpha.
\end{eqnarray}
 Taking another covariant
derivative at $p$, it gives
\begin{eqnarray}\label{nabla2y}
  \nabla_{e_k}^\perp\nabla_{e_i}^\perp y^\perp&=&-e_k\langle y, e_j\rangle h_{ij}^\alpha e_\alpha-\langle y, e_j\rangle h_{ij, k}^\alpha e_\alpha\nonumber\\
  &=&-\langle y, h_{kj}^\beta e_\beta\rangle h^\alpha_{ij}e_\alpha-\langle y, e_j\rangle h^\alpha_{ik, j}e_\alpha,
\end{eqnarray}
where we used (\ref{nablay}), $\nabla_{e_k}e_j(p)=0$, and the
Codazzi equation in the last equality. Taking the trace of
(\ref{nabla2y}) and using $H=T^\perp$, we conclude that
\begin{eqnarray*}
  \Delta^\perp y^\perp&=&-\langle y, h_{ij}^\beta e_\beta\rangle h^\alpha_{ij}e_\alpha-\langle y, e_j\rangle H^\alpha_{,j} e_\alpha\\
  &=&-y^\beta h_{ij}^\beta h_{ij}^\alpha e_\alpha+\langle y, e_j\rangle\langle T, e_i\rangle h_{ij}^\alpha e_\alpha\\
  &=&-y^\beta h_{ij}^\beta h_{ij}^\alpha e_\alpha-\langle T, e_i\rangle
  \nabla_{e_i}^\perp y^\perp\\
  &=&-y^\beta h_{ij}^\beta h_{ij}^\alpha e_\alpha-\nabla_{T^T}^\perp y^\perp.
\end{eqnarray*}
This proves (\ref{Ly}).

Since $\Delta x=H$ and $H=T^\perp$, we have
\begin{eqnarray*}
  \Delta x^A=\langle H, E_A\rangle=\langle T^\perp, E_A\rangle=\langle T,
  E_A^\perp\rangle=\langle T, E_A\rangle-\langle T, E_A^T\rangle=T^A-\langle T, E_A^T\rangle.
\end{eqnarray*}
Hence
\begin{eqnarray*}
  \mathcal{L}x^A=\Delta x^A+\langle T, \nabla x^A\rangle=\Delta
  x^A+\langle T, (E_A)^T\rangle=T^A,
\end{eqnarray*}
This proves (\ref{LxA}). \hfill Q.E.D.

\subsection{Hamiltonian L-stability}
In this subsection, we will define Hamiltonian L-stability of
Lagrangian translating solitons. Recall the definition of
Hamiltonian variations on a Lagrangian submanifold.
\begin{definition}\cite{O}
  Let $(M, \bar\omega)$ be a symplectic manifold $M$. Let $\Sigma\subset M$ be a Lagrangian submanifold and $V$ be a vector field along $\Sigma$. $V$
  is called a Hamiltonian variation if it satisfies that the one form $i^*(V\rfloor\bar\omega)$ on $\Sigma$ is exact.
\end{definition}
The Hamiltonian variation has an equivalent definition.
\begin{lemma}\label{Hv}\cite{O}
  A normal variation $V$ on $\Sigma$ is Hamiltonian if and only if
  \begin{eqnarray*}
    V=J\nabla f,
  \end{eqnarray*}
  where $f$ is a function on $\Sigma$ and $\nabla$ is the gradient on $\Sigma$ with respect to the induced metric.
\end{lemma}

Now we are ready to define Hamiltonian L-stability of Lagrangian
translating solitons.
\begin{definition}\label{LHs}
  We say a Lagrangian translating soliton $\Sigma$ is Hamiltonian L-stable if for every compactly supported Hamiltonian
  variations $\Sigma_s$ with $\Sigma_0=\Sigma$, $F''=\int_\Sigma-\langle V, LV\rangle e^{\langle T, x\rangle}\geq 0$.
\end{definition}

\section{Proof of the main theorem}

Note that the normal bundle brings much difficulty to the study of
L-stability of translating solitons in the general
higher codimension case. However, in \cite{O}, Oh studied
Hamiltonian stability of minimal Lagrangian submanifolds in
K\"ahler-Einstein manifolds, and characterized Hamiltonian stability
by a condition on the first eigenvalue of $\Delta$ acting on
functions. The key point of Oh's proof is that, for a minimal
Lagrangian submanifold of a K\"ahler-Einstein manifold, the set of
Hamiltonian variations is an invariant subspace of the Jacobi
operator. This idea was then used to study Hamiltonian
F-stability of Lagrangian self-shrinkers (\cite{LZ}, \cite{Y}). It
is natural to think that this propery also holds for Lagrangian
translating solitons. This property inspired us to show the
following equality, which well characterizes how the operator $L$
acts on Hamiltonian variations.

\begin{theorem} Suppose $\Sigma^n\subset\textbf{C}^n$ is a Lagrangian translating soliton. Then for every function $f$ on $\Sigma$,
  \begin{eqnarray}\label{e3.1}
    LJ\nabla f=J\nabla {\mathcal{L}} f.
  \end{eqnarray}
  This implies that the set of Hamiltonian variations is an invariant subspace of the operator $L$.
\end{theorem}

{\it Proof. } Fix a point $p$. We choose a local orthonormal basis $\{e_i\}_{i=1}^n$ of $T\Sigma$ such that $\nabla_{e_i}e_j(p)=0$. Then since $\Sigma$ is Lagrangian, $\{e_{n+i}=Je_i\}_{i=1}^n$ is a local orthonomal basis of $N\Sigma$. In the following we compute at the point $p$. It is easy to compute that
\begin{eqnarray*}
LJ\nabla f&=&\Delta^\perp(J\nabla f)+\nabla^\perp_{T^T}(J\nabla f)+h^{n+k}_{il}f_kh_{il}^{n+j}Je_j\nonumber\\
&=&\left(f_{jii}+\langle T, e_k\rangle f_{jk}+f_kh^{n+l}_{ik}h^{n+l}_{ij}\right)Je_j,
\end{eqnarray*}
where in the last equality we used the Lagrangian property $h_{il}^{n+k}=h_{ik}^{n+l}$.
On the other hand,
\begin{eqnarray*}
  J\nabla{\mathcal{L}}f&=&J\nabla\left(\Delta f+T^T f\right)\nonumber\\
  &=&f_{iij}Je_j+J\nabla\langle T^T,\nabla f\rangle\nonumber\\
  &=&\left(f_{iji}-f_iR_{jkik}+e_j\langle T^T, \nabla f\rangle\right)Je_j\nonumber\\
  &=&\left(f_{jii}-f_ih_{ij}^{n+l}h_{kk}^{n+l}+f_i h_{jk}^{n+l}h_{ik}^{n+l}+\langle\nabla_{e_j}T^T, \nabla f\rangle+\langle
  T^T, \nabla_{e_j}\nabla f\rangle\right)Je_j\nonumber\\
  &=&\left(f_{jii}-f_ih_{ij}^{n+l}\langle T^\perp, e_{n+l}\rangle+f_kh_{ik}^{n+l}h_{ij}^{n+l}+\langle\overline{\nabla}_{e_j}T,\nabla f\rangle-\langle\overline{\nabla}_{e_j}T^\perp, \nabla f\rangle\right.\nonumber\\
  &&\left.+f_{jk}\langle T, e_k\rangle\right) Je_j\nonumber\\
  &=&\left(f_{jii}-f_ih_{ij}^{n+l}\langle T^\perp, e_{n+l}\rangle+f_kh_{ik}^{n+l}h_{ij}^{n+l}+\left\langle T^\perp, \overline{\nabla}_{e_j}(f_ke_k)\right\rangle+\langle T, e_k\rangle f_{jk}\right) Je_j\nonumber\\
  &=&\left(f_{jii}-f_ih_{ij}^{n+l}\langle T^\perp, e_{n+l}\rangle+f_kh_{ik}^{n+l}h_{ij}^{n+l}+f_k\langle T^\perp, h_{jk}^{n+l}e_{n+l}\rangle+\langle T, e_k\rangle f_{jk}\right)Je_j\nonumber\\
  &=&\left(f_{jii}+\langle T, e_k\rangle f_{jk}+f_k h_{ik}^{n+l}h_{ij}^{n+l}\right)
  Je_j,
\end{eqnarray*}
where in the third equality we used the Ricci formula; in the fourth
equality we used the Gauss equation; and in the fifth equality we
used the translating soliton equation $H=T^\perp$. This proves the
theorem. \hfill Q.E.D.

Now we recall our main theorem.
\begin{theorem}
Any Lagrangian translating soliton is Hamiltonian L-stable.
\end{theorem}
{\it Proof.} Recall that the second variation formula for $F$ is
\begin{eqnarray}\label{e4.1}
  F''=\int_\Sigma-\langle V, LV\rangle e^{\langle T, x\rangle}.
\end{eqnarray}
Now Assume $V$ is a compactly supported Hamiltonian variation, then
there exists a function $f$, such that $V=J\nabla f$. Putting it
into (\ref{e4.1}), and using (\ref{e3.1}) , we have
\begin{eqnarray}\label{e4.2}
  F''=\int_\Sigma-\langle J\nabla f, LJ\nabla f\rangle e^{\langle T, x\rangle}=\int_\Sigma-\langle J\nabla f, J\nabla\mathcal{L}f \rangle e^{\langle T, x\rangle}
  =\int_\Sigma(\mathcal{L}f)^2 e^{\langle T, x\rangle}\geq 0,
\end{eqnarray}
where the last equality used Lemma \ref{uLv} and the fact that
$V=J\nabla f$ is compactly supported implies that $\mathcal{L}f$ is
compactly supported. This proves the theorem. \hfill Q.E.D.

\end{document}